\documentclass[a4paper,11pt]{article}
\usepackage[a4paper,bindingoffset=0.2in,
            left=0.8in,right=0.9in,top=1in,bottom=1in,
            footskip=.25in]{geometry}
\usepackage{amsmath}
\usepackage{amsfonts}
\usepackage{amssymb}
\usepackage{amsthm}

\begin{document}

\begin{center}
{\LARGE  THE CDE METHOD}
\end{center}

\begin{center}
A TECHNIQUE IN FUNCTIONAL EQUATIONS \\
\end{center}

\begin{center}
\ \\
Athanasios Kontogeorgis and Rafail Tsiamis \\
\end{center}

\begin{center}
January 30th, 2019
\end{center}

\begin{center}
\ \\
\ \\
\ \\
\ \\
\ \\
\ \\
\textbf{Abstract}

\ \\

In this article we present an extremely effective and relatively unknown approach to solving functional equations that appear in mathematical competitions. We aim to explain the philosophy of this novel method through numerous examples, which also highlight how this idea can be paired with other useful techniques to crack challenging problems. 
\end{center}

\newpage 

\textbf{1. Introduction} \\

The method \textit{cde} is named after the main relation it aims to obtain: 
\begin{center}
$f(x+c) = f(x)+d,$ for all $x \geq e.$ 
\end{center}
Next, we consider some difference of the form $P(x+c,y) - P(x,y)$ (or $P(x,y+c) - P(x,y),$ or something similar) in order to create ''coefficient imbalances'', thus getting things to cancel out, which leads to useful relations. 

The idea of creating a coefficient imbalance truly lies at the core of our method; sometimes it comes up in exercises that do not explicitly use the cde technique in similar ways - for example, getting $f(cx) = df(x),$ for all $x\geq e,$ and then considering a difference of the form $P(cx,y) - dP(x,y)$ or a similar quotient. The good thing about this method is that it's routine, aimless (we don't really know what we're going to get, but it's more often than not something useful) and very effective, even up against challenging problems.

This concept is used in conjunction with another very strong idea: ''Iterating'' the given relation. When the condition states $f(A(x,y)) = B(x,y)$ and the argument $A(x,y)$ includes $f(x)$ (respectively, $f(y)$), then it is often useful to make the substitution $x\to A(x,z)$ (respectively, $y\to A(y,z)$).

In general, a relation of the form $f\left(f(x) + S_1\right) = x + S_2$  is always going to give a cde-relation, simply by iterating $x\to f(x)+S_1$ to obtain $f(x+S_1+S_2) = f(x)+S_1+S_2.$ \\

\textbf{Iteration lemma}: $f(f(x)+S_1) = x+S_2$ implies $f(x+K) = f(x) + K,$ where $K = S_1+S_2.$ It is trivial to see inductively that $f(x+c) = f(x) + d$ implies $f(x+Nc) = f(x) + Nd,$ for any natural number $N.$ 
\\
\\ A simple example of this idea: 
Determine all functions $f: \Bbb{Q}\to \Bbb{Q},$ such that 
\begin{center}
$f(f(x)+y) = x+f(y).$
\end{center}
Iterate by taking $x\to f(x)+z,$ which gives $f(x+y+f(z)) = f(x) + f(y) + z$ (1). However, the initial equality also gives $f(f(z)+(x+y)) = z + f(x+y)$ (2). From (1),(2), it follows that $f(x+y) = f(x)+f(y),$ which is Cauchy 's equation, well-known to have the solution $f(x)=cx$ over the rationals. Substituting back to the initial equation, we find $c=1$ or $c=-1.$ \\

Another idea this method relates very well to is proving that a function is either injective or periodic; we assume that $f(a)=f(b)$ for some $a,b$ and obtain a relation of the form $f(x) = f(x+T),$ where $T=T(a,b).$ Substituting $x\to x+T$ in the initial relation usually either gives that $T=0$ (which forces certain conditions on $a,b$) or reduces to only trivial solutions.

Lastly, an idea that works very well independently of our technique: trying to get some arguments between the two sides equal, which would lead to terms cancelling out. In particular, when the range of the function is the positive reals, some terms can never become zero, which imposes certain inequalities on the values of $f.$

\newpage

\textbf{2.	Example problems} \\
\\ Unless otherwise stated, we denote the given relation for numbers $x,y$ as $P(x,y).$ We also assume that the given relation holds for all $x,y$ in the domain of the function. 
\\ 
\\
\\ \textbf{1.}	Determine all functions $f: \Bbb{R}\to \Bbb{R}$ such that 
\begin{center}
 $f(f(x)+xy) + f(y^2) = 2017x + yf(x+y).$  \\~\\
\end{center} 
   \textit{Solution:} $P(0,0): f(f(0)) + f(0) = 0$ \hfill (1)
\\ $P(0,y): f(f(0)) + f(y^2) = yf(y),$ thus $f(y^2) = yf(y) + f(0).$ \hfill (2)
\\ Setting $y=1,$ we get that $f(0) = f(f(0)) = 0.$
\\ $P(x,0): f(f(x)) = 2017x$.  Also, $f(2017x) = f(f(f(x))) = 2017f(x).$ \hfill (3) 
\\ Furthermore, if $f(u)=0,$ for some $u,$ it follows that $2017u = f(0) = 0,$ hence $u=0.$ 
\\ Therefore $f(u) = 0 \iff u=0.$ 
\\ Now we create the desired 'coefficient imbalance' by $P(2017x, 2017y)$:
\\ $f(f(x)+2017xy) + 2017f(y^2) = 2017x + 2017yf(x+y).$ 
\\ Equivalently, $f(f(x)+2017xy) - 2017x = 2017(yf(x+y) - f(y^2)) = 2017(f(f(x)+xy) - 2017x),$ using $P(x,y).$ 
\\ That is, using (3), $f(f(x)+2017xy) = 2017f(f(x)+xy) - 2016f(f(x)).$ \hfill (4)   
\\ By (3), we see that $f(x)$ spans all nonzero real numbers for nonzero $x$. Thus, we may set $f(x)=z$ where $z$ is any nonzero real number. Also, set $\displaystyle y \to \frac{y}{x}$ ($x$ is nonzero). 
\\ $(4)$ becomes: $f(z+2017y) = 2017f(y+z) - 2016f(z),$ that is \\
$Q(x,y): f(x+2017y) + 2016f(x) = 2017f(x+y)$ for any real $y$ and any nonzero $x.$ By (3), this also holds for $x=0.$
\\ $Q(2017x,y): f(x+y) + 2016f(x) = f(2017x+y).$
\\ Hence $Q(y,x) + Q(2017x,y)$ gives $f(x)+f(y) = f(x+y),$ for all $x,y$ in $\Bbb{R}.$ \hfill (5)
\\ Using (2),(3),(5), $P(x,y)$ now gives $f(xy) = yf(x),$ which implies $xf(y) = f(xy) = yf(x),$ or $f(x)=cx.$ By (3), we have $c^2 = 2017.$ 
\\ Thus the two solutions are $f(x) = \sqrt{2017} x$ and $f(x)=-\sqrt{2017}x.$  \\~\\~\\
\textbf{2.}	Determine all functions $f:\Bbb{R}^+\to \Bbb{R}^+$ such that 
\begin{center}
$f(x+f(x)+y) = f(2x)+f(y).$ \\~\\
\end{center}
\textit{Solution:} We try to make arguments equal: $x+f(x)+y = 2x$ would require $y = x-f(x).$ Hence, assume that there is some $x_0$ such that $x_0 > f(x_0).$ Then $P(x_0,x_0-f(x_0))$ gives $f(x_0-f(x_0))=0,$ clearly absurd. This means that $f(x) \geq x.$ 
\\ On the other hand, we need only prove that f is injective. Indeed, $P(x,2y)$ and $P(y,2x)$ would then imply that $f(x+f(x)+2y) = f(2x) + f(2y) = f(y+f(y)+2x),$ thus $f(x) = x+c$ by injectivity; we verify that $f$ satisfies the initial relation for any nonnegative constant $c.$ 
\\It suffices to prove that $f$ is injective. Assume for contradiction that there exist $a>b$ such that $f(a)=f(b)$; then $P(a,b), P(b,a)$ give $f(2a) = f(2b).$ Now $P(a,x)$ and $P(b,x)$ give \\ $f(a+f(a)+x) = f(b+f(b)+x),$ or $f(x) = f(x+T)$ for $T= a-b > 0, x > b+f(b).$ Inductively, this gives $f(x) = f(x+nT).$
However, for sufficiently large $n,$ such that $x+nT > f(x),$ we would have $f(x+nT) \geq x+nT > f(x) = f(x+nT),$ which is absurd; it follows that $a=b,$ implying that the function is injective, as desired. \\~\\~\\
\textbf{3.}	(Iran TST) \ Determine all functions $f: \Bbb{N}^* \to \Bbb{N}^*$ such that 
\begin{center}
$(f(a)+b)f(a+f(b)) = (a+f(b))^2.$ \\~\\
\end{center}
\textit{Solution 1:} First of all, observe that the function is injective; indeed, the initial relation gives $f(a)+b \mid (a+f(b))^2.$ If $f(a_1) = f(a_2) = c,$ it follows that $x+c$ divides \\
$\gcd((a_1+f(x))^2,(a_2+f(x))^2) = \gcd(a_1+f(x),a_2+f(x))^2 = \gcd(|a_1-a_2|,a_2+f(x))^2$; \\ in particular, it divides $(a_1-a_2)^2,$ which is absurd for sufficiently large $x,$ unless $a_1=a_2.$
\\ Now consider $P(p-f(x), x)$ for some prime number $p$ and $x>1$; it gives $(f(p-f(x))+x)f(p) = p^2,$ therefore $f(p)=1$ or $f(p)=p$ (since $f(p-f(x)) + x >1$). Since the function is injective, $f(p)=1$ for at most one number, thus we must have $f(p)=p$ for all sufficiently large primes $(p>M),$ which also gives $f(p-f(x))+x = p,$ or $f(p-f(x)) = p-x$ for any $x$ and any prime $p>M$ (1).
\\Now we iterate (1) by taking $x\to q-f(x),$ where $q$ is a prime such that $p > q >M.$ This gives:
\\$f(x+(p-q)) = f(x) + (p-q)$ for any sufficiently large $x>N,$ hence $f(x+K) = f(x)+K$ for any $x>N.$ 
\\Since f is injective and defined on the natural numbers, it follows that there exists a constant $S,$ such that $f(x) > N$ for any $x>S.$ Taking $a,b > \max\{S,N\}$ gives $\min\{a, f(a), b, f(b)\} > N.$ Hence, considering $P(a+K,b) - P(a,b),$ we obtain $(f(a)+b) + f(a+f(b)) = 2(a+f(b))$; along with $P(a,b),$ this gives $\{f(a)+b, f(a+f(b)\} = \{a+f(b)\},$ thus $f(a) = a+c.$ \\ Substituting back to $P(a,b),$ we can see that $f(n) = n,$ for all natural numbers, is the only solution.\\ \\
\textit{Solution 2:} $P(a,a)$ gives $f(a+f(a)) = a+f(a)$ (1) for all natural numbers $a.$ \\ Denote by $F$ the set of fixed points of $f: F = \{n \in \Bbb{N}^* : f(n)=n\}.$ By (1), $a+f(a) \in F$ for all $a\in \Bbb{N}^*.$ Take a fixed point $A$; by $P(A,A)$, $2A\in F.$ Now:
\\$P(x, 2A): (f(x)+2A)f(x+2A) = (x+2A)^2$
\\$P(x+A, A): (f(x+A)+A)f(x+2A) = (x+2A)^2.$
\\It follows that $f(x+A) = f(x) + A$ for all $x\in \Bbb{N}^*, A \in F.$
\\$P(x, A): (f(x)+A)^2 = (f(x)+A)f(x+A) = (x+A)^2,$ hence $(f(x)-x)(f(x)+x+2A) = 0.$ It follows that $f(x) = x$ for all positive integers $x.$  \\~\\~\\
\textbf{4.} (Iran Training Camp)	Determine all functions $f: \Bbb{N}^*\to \Bbb{N}^*,$ such that 
\begin{center}
$f(m^2+f(n)) = f(m)^2+n.$ \\~\\
\end{center}
\textit{Solution:} By the iteration lemma applied to $P(1,n),$ we get $f(n+A) = f(n) + A,$ where $A = f(1)^2 + 1.$ 
Now $P(m+A, n) - P(m,n): f(m^2 + 2Am + A^2 + f(n)) - f(m^2+f(n)) = 2Af(m) + A^2.$\\
Now we can extract the terms $2Am + A^2$ from the arguments, since they are natural multiples of $A.$ But, when we cancel the $f(m^2+f(n)),$ this only leaves $f(m) = m,$ for all $m.$ \\
Seriously, this should feel like cheating. \\~\\~\\
In the following problems, we may not look for a relation of the form $f(x+c) = f(x)+d$; however, it is important to understand the \textbf{"cde philosophy"}: we substitute constants to break terms into simpler ones, for example of the form $f(x^2 + \text{something}) = xf(x) + \text{something else}$, etc.  The main goal is to get a Cauchy-type equation plus constants, which we use to reduce the initial relation. Another problem we may encounter is neatly having $x$ (respectively $y$) in the one argument, yet having a complex function of them on the other side. This problem is tackled by subtracting $P(x,a) - P(x,b)$ (respectively $P(a,y) - P(b,y)$) and treating the resulting functions of $a,b$ as $c$ and $d.$ \\~\\
These underlying ideas are best understood through some examples: \\~\\~\\
\textbf{5.} (Italian TST for IMO 2016)	Prove that there does not exist a function $f:\Bbb{R}^+ \to \Bbb{R}^+$ such that 
\begin{center}
$f(f(x)+y) = f(x)+3x+yf(y).$ \\~\\
\end{center}
\textit{Solution:} First of all we investigate what equal arguments lead to: if, for some $x_0,$ $x_0 > f(x_0),$ then $P(x_0, x_0- f(x_0))$ gives $3x_0 + (x_0-f(x_0))f(x_0-f(x_0)) = 0,$ clearly absurd. Thus $f(x)\geq x,$ for all $x.$ \\ 
Now, we observe that the $y$ in the left-hand argument is a great candidate for our method, however the term $yf(y)$ on the right-hand side is cumbersome. As suggested before, this prompts us to take the difference $P(s,y)- P(t,y)$ for some $s,t$ to be specified later.
\\ We obtain $f(f(s)+y) - f(f(t)+y) = f(s) - f(t) + 3(s-t).$ Hence, if we denote \\ $c = s-t, d = f(s) - f(t),$ this gives $f(y+d) = f(y) + 3c + d$ for all $y > f(t).$ To avoid complications, we can also assume that $s > 2f(t),$ so that $c = s-t > f(t) + t - t > 0,$ and $d = f(s) - f(t) \geq  s – f(t) > f(t) > 0.$ \\ 
Now subtract $P(x,y+d) - P(x,y)$. This gives $d+3c = df(y) + (d+3c)y + d(d+3c)$ \hfill (1) \\
By our previous observations, $c,d > 0,$ thus (1) clearly fails for sufficiently large $y;$ we are done. \\~\\~\\
\textbf{6.} Determine all functions $f: \Bbb{R}^+ \to \Bbb{R}^+$ such that 
\begin{center}
$f(xf(y)+f(x)) = 2f(x)+xy.$ \\~\\
\end{center}
\textit{Solution:} If we let $f(1)=a,$ then $P(1,x)$ gives $f(f(x)+a) = x+2a,$ thus $f$ assumes all values on $(2a, +\infty).$ Also, the iteration lemma gives $f(x+3ka) = f(x)+3ka$ \hfill (*). \\ Now, the familiar step $P(x,y+3ka) - P(x,y)$ gives $f(xf(y)+f(x)+3kax) - f(xf(y)+f(x)) = 3kax.$  \\ 
Now, let $x\to  \frac{x}{3ka}$ be any positive real number. Then, since $f(y)$ may assume any value greater than $2a,$ we see that $\displaystyle \frac{xf(y)}{3k}+f\left(\frac{x}{3k}\right)$ may assume any value $z$ greater than $\displaystyle g(x) = \frac{2ax}{3k} + f\left(\frac{x}{3k}\right).$ \\ Hence, $f(z+x) = f(z)+x$ for all $x\in \Bbb{R}^+, z > g(x).$ But if this holds for $(z,x),$ it also holds for $(z,x+3ka),$ by (*).\\ Taking $k$ to be sufficiently large, we actually get $f(z+x) = f(z)+x$ for all $x,y.$ \\ Thus $f(x)+z = f(x+z) = f(z)+x$ for all $x,y,$ which gives $f(x)=x+c.$ Substituting this into the initial equation gives $c=1,$ i.e. $f(x)=x+1$ is the unique solution to this equation. \\~\\~\\
\textbf{7.}	Determine all functions $f: \Bbb{R}^+\to \Bbb{R}^+$ such that 
\begin{center}
$f(xf(y)+x) = xy + f(x).$ \\~\\
\end{center}
\textit{Solution:} By fixing $x$ and letting $y$ span $\Bbb{R}^+,$ we obtain that $f$ is surjective onto an interval of the form $(M, +\infty),$ where $M = f(x).$ Taking $x\to1,$ we obtain $f(f(y)+1) = y + f(1),$ thus $f(y+c) = f(y) + c$ where $c = 1+f(1)$ by the iteration lemma. \\ 
$P(x, y+c) - P(x,y)$ gives $f(xf(y)+x+cx) - f(xf(y)+x) = cx.$ Since $f(y)$ can be equal to any sufficiently large positive real number, it follows that $xf(y)+x$ can also be equal to any sufficiently large $z.$ \\ Now, taking $\displaystyle x\to \frac{x}{c}, \ xf(y)+x \to z,$ we obtain $f(x+z) = f(x) + z$ for all $x \in \Bbb{R}^+, z > M(x),$ where $M(x)$ is a function of $x.$ \\ 
Take $z > \max \{M(x),M(y)\}.$ Then $f(x+y+z) = f(x) + y+z$ (since $y+z > M(x)$), but also
$f(x+y+z) = f(y) + x+z.$ It follows that $f(x) = x+a,$ for all $x$ and some constant $a.$ \\ Substituting back to the original equation, this gives $a=0,$ i.e. $f(x) = x.$ \\ \\
Otherwise, we could have obtained this by seeing that, since $f(x+kc) = f(x)+kc,$ \\ $f(x+kc+z) = f(x)+z+kc = f(x+kc)+z,$ and $kc$ may be made sufficiently large, so that in fact $f(x+z) = f(x)+z,$ for all $x,z  \in \Bbb{R}^+,$ thus again $f(x)=x+a$ leading to $f(x)=x.$ \\~\\~\\
\textbf{8.}	Determine all functions $f: \Bbb{R}^+ \to \Bbb{R}^+$ such that 
\begin{center}
$f(x^2+f(y)) = xf(x) + y.$ \\~\\
\end{center}
This can be thought of as the model application of the cde method: first we prove that some values span all sufficiently large numbers, and we isolate them; we substitute the relations back to the initial equation to get a Cauchy-type equation; step by step we reduce it to Cauchy + a constant which eventually gives the solution for all sufficiently large numbers; we then conclude for smaller values. We really don't care about constants inside arguments or ''sufficiently large'' things because they can be eventually reduced to simpler expressions. \\ \\
\textit{Solution:} First of all, by fixing $x$ and varying $y,$ we obtain that $f$ is surjective onto an interval of the form $[M, +\infty)$ due to the right-hand side.\\
$P(x,a)$ gives $f(x^2 + b) = xf(x) + a$ where $b=f(a)$ \\
$P(u,y)$ gives $f(f(y) + c) = y + d,$ where $c=u^2, d=uf(u)$ \\
Thus $P(x,y)$ becomes $f(x^2+f(y)) = [f(x^2+b) - a] + [f(f(y)+c) - d].$ \\ By setting $x\to \sqrt{x}, \ f(y)\to z$ (for any $z\geq M$, by the surjectivity proved above) and $a+d = D,$ we obtain: \\
$f(x+y) = f(x+b) + f(y+c) - D$ for any $x  \in \Bbb{R}^+, y>M.$ \\~\\ 
Now let's prove a very useful lemma: \\ \\ 
\textbf{Cauchy-type lemma (CTL)}: Suppose that $f(x+y) = f(x+u) + f(y+v) + w$  \\ for constants $u,v,w$ and all sufficiently large $x,y.$ Then $f(x) = c(x) + d,$ for all sufficiently large $x$ and some constant $d,$ where $c(x)$ is a solution to Cauchy's equation $c(x+y) = c(x) + c(y).$\\  \\ 
\textit{Proof:} Everything we do assumes sufficiently large $x,y.$ Denote the given relation by $A(x,y).$
If $u=v,$ proceed to the next step. Otherwise, $A(x,y) - A(y,x)$ gives $f(x+u) = f(x+v) + s$ for some $s,$ hence $A(x,y)$ becomes $f(x+y) = f(x+v) + f(y+v) + w',$ where $w' = w+s.$ \\
Now $A(x+v,y) - A(x,y+v)$ gives $f(x+2v) = f(x+v) + t$ for some $t,$ implying $f(x+v) = f(x) + t$ for all sufficiently large $x.$ \\
Now $A(x,y)$ becomes $f(x+y) = f(x) + f(y) + w'',$ where $w'' = w'+s+2t.$ Taking $c(x) = f(x) + w'',$ this gives $c(x+y) = c(x) + c(y$) for all sufficiently large $x,y.$ The desired lemma is proved. \\~\\
We return to our problem.\\ By the lemma, $f(x) = c(x) + A$ for all sufficiently large $x$ and some constant $A.$ \\  Since $f$ maps to positive reals, $c$ is bounded below; it is well known that $c(x)$ being Cauchy and lower-bounded means $c(x) = cx$ for all sufficiently large $x.$ 
\\This means that $f(x) = cx+A$ for all sufficiently large $x,$ say $x>M.$
\\Taking $x = M+1,$ we also have $x^2 > M,$ thus $f(x^2+f(y)) = c(x^2+f(y)) + A.$
\\Therefore $cx^2 + cf(y) + A = cx^2 + Ax + y.$ Taking $x$ to infinity, this means that $A=0.$ \\ Hence, $\displaystyle f(y) = \frac{y}{c},$ for all $y,$ in particular $cy = f(y) = \displaystyle  \frac{y}{c},$ for sufficiently large $y$ gives $c=1$ (since $f$ is positive), i.e. $f(x) = x,$ for all sufficiently large $x.$ \\ Now $P(M+1,y)$ implies $f(y)=y,$ for all $y.$ \\~\\~\\
\textbf{9.} (Turkey TST 2014)	Determine all functions $f: \Bbb{R}\to \Bbb{R}$  such that 
\begin{center}
$f(f(y) + x^2 + 1) + 2x = y + f^2(x+1).$ \\~\\
\end{center}
\textit{Solution:} Clearly, allowing $y$ to span all of $\Bbb{R}$ on the RHS gives that $f$ is surjective. \\  Furthermore, taking $P(x, y_1), P(x, y_2)$ for $y_1, y_2$ such that $f(y_1) = f(y_2)$ yields $y_1 = y_2,$ meaning that the function is also injective. \\ 
$P(0,y)$ gives $f(f(y)+1) = y + f^2(1),$ thus $f(x+c) = f(x)+c$ by the iteration lemma, where $c = 1+f^2(1)$ is positive. \\ \\
We note that it is possible to find an elegant solution by iterating $y\to f(y)+1$ back to the initial equation; we leave this to the interested reader and showcase a routine, brute-force solution using the cde method instead. \\ \\
We will mostly be spamming $x\to x+c$ in relations in order to create 'coefficient imbalances' that allow us to isolate things after considering $P(x+c, y) - P(x,y).$ \\
$P(x+c, y) - P(x,y):  f(x^2+f(y)+1 + 2cx + c^2) + 2c - f(x^2+f(y)+1) = c^2 + 2cf(x+1).$ \\  
But since $f(y)$ is surjective onto $\Bbb{R}$, $x^2+f(y)+1$ can also be any real value $z$ (by varying $y$ for a fixed $x$). Hence, since $f(A)+2c = f(A+2c),$ we obtain $f(c^2+2c(x+1)+z) = c^2 + f(z) + 2cf(x+1).$ Taking $x\to  x-1, z\to y,$ this can be written as $Q(x,y):$ \\ 
$Q(x,y): f(c^2 + 2cx + y) = f(y) + 2cf(x) + c^2.$ \\
$Q(0,y): f(c^2+y) = f(y) + c^2 + 2cf(0)$ \hfill (1) \\ Now $Q(x,y)$ becomes: $f(2cx+y) + 2cf(0) = f(y) + 2cf(x).$ \\ 
In particular, $Q(x+c,0) - Q(x,0)$ yields $f(2cx+2c^2) - f(2cx) = 2c^2.$\\ By (1), the LHS equals $2c^2 + 4cf(0).$ Since $c>0,$ it follows that $f(0) = 0.$ \\ Finally, $Q(x,y)$ can be written as $R(x,y):$
$f(2cx+y) = f(y) + 2cf(x).$ \\ For $y=0,$ we get $2cf(x) = f(2cx),$ thus $f(x+y) = f(x)+f(y).$ 
\\We use this in the initial equation: $f(f(y))+f(x^2)+2x = y + f(x)^2 + 2f(x);$ let this be $S(x,y).$ \\ 
$S(0,y)$ immediately gives $f(f(y)) = y$ for all $y.$  \\ In all, our function has the following properties:\\ 
(i) $f(x+y) = f(x) + f(y)$, \\ (ii) $f(f(x)) = x$, \\ (iii) $f(x^2) + 2x = f(x)^2 + 2f(x)$. \\ The first property gives $f(x) = ax$ for all rationals (well-known). Using this in (iii), we must have $a=1.$
We create our beloved coefficient imbalance to finish: (iii) for $2x$ gives \\ $4f(x^2) + 4x = 4f(x)^2 + 4f(x),$ thus $f(x^2) + x = f(x)^2 + f(x).$ Combined with (iii), this immediately gives $f(x) = x.$\\~\\~\\
\textbf{10.} Determine all functions $f: \Bbb{R}\to\Bbb{R}$ such that 
\begin{center}
$f(x^n +f(y)) = f(x)^n + y,$
\end{center}
for some fixed positive integer $n \geq 2.$\\~\\ 
\textit{Solution:} Fix $x$ and let $y$ span all of $\Bbb{R}$; this proves that $f$ is surjective. Also, if $f(y_1) = f(y_2)$, it is immediate that $y_1 = y_2$, thus f is injective. Let $f(0)= a$, and let $b$ be such that $f(b) = 0.$ 
\\ $P(x,b): f(x^n) = f(x)^n + b$
\\ $P(0,y): f(f(y)) = y + a^n$
\\ Now $P(x,y)$ becomes $f(x^n + f(f(y)) = f(x^n) + f(f(y)) - (a^n+b)$. Writing $a^n + b = D$, \\ $f(y)\to y \in \Bbb{R}$ and $x^n \to x  \in \Bbb{R}^+ \cup \{0\},$ we get:
\\ $Q(x,y): f(x+y) = f(x) + f(y) - D$ for all nonnegative $x$ and real $y$. \\ Let $g(x) = f(x) - D,$ then
$Q(x,y): g(x+y) = g(x) + g(y)$ for all nonnegative $x$ and real $y.$ $Q(x,0)$ gives $g(0) = 0,$ hence $Q(x,-x)$ implies $g(-x) = -g(x)$. Thus $g(x+y) = -g((-x)+(-y)) = g(x) + g(y)$ also holds for $x,y$ both negative; it follows that $f(x) = g(x) + D$ for all real $x,$ where $g$ is additive.
\\ $P(x,y)$ can now be written as $g(x^n) + g(g(y)) + g(D) + D = (g(x) + D)^n + y.$ Since $g(0) = 0$, $P(0,y)$ gives $g(g(y)) + g(D) + D = D^n + y$ \hfill (1).
\\$P(x,y)$ becomes $g(x^n) = (g(x)+D)^n - D^n$ \hfill (2).
\\Since $g(q) = aq$ for all rational numbers $q$, (1) gives $a^2 = 1$, that is $a=1$ or $a=-1.$ 
\\Applying (2) for rational $x$, we get $\displaystyle ax^n = \left(\sum_{k=0}^n \binom{n}{k} (ax)^k D^{n-k}\right)-D^n$,  which implies \\ $\displaystyle \left(\sum_{k=1}^{n-1} \binom{n}{k}(ax)^k D^{n-k}\right)+x^n (a^n-a)=0.$
\\Since this is a polynomial equality that holds for infinitely many $x$ (all the rational numbers), it must hold identically; it follows that all the coefficients of this polynomial must be identically zero.
\\Observe that each coefficient may only be $0$ if and only if $D=0.$ Thus $f(x)$ is additive.
\\We furthermore obtain the relations $f(f(x)) = x$ and $f(x^n) = f(x)^n$ and $f(q) = q$ or $f(q) = -q$ on the rationals. For even $n,$ $f(q) = -q$ doesn't work; for odd $n,$ $f$ being a solution implies that $-f$ is also a solution, since $f(-x) = -f(x)$ because $f$ is additive. Hence, we may assume that $f(q) = q$ for all rational numbers.
\\Now $f(x+q) = f(x) + f(q) = f(x) + q.$ Thus $f(x^n) = f(x)^n$ gives $f((x+q)^n) = [f(x)+q]^n,$ therefore:
$\displaystyle \sum_{k=0}^n \binom{n}{k}q^kf(x^{n-k})=\sum_{k=0}^n\binom{n}{k}q^kf(x)^{n-k}$,\\ by repeatedly using the additivity of $f$ to extract terms from the brackets. Fixing $x,$ we may interpret this as a polynomial equality in $q$ that holds for infinitely many $q,$ (all the rational numbers) thus identically; it follows that the coefficients of any degree must be equal, meaning $f(x^k) = f(x)^k$ for all $k\leq n.$ In particular, $f(x^2) = f(x)^2\geq 0$ proves that $f$ is nonnegative for all nonnegative numbers; since it is also additive, it is well-known that $f(x) = cx$ everywhere. \\Verifying back to the initial equation, it follows that $f(x) = x$ if $n$ is even, and if $n$ is odd, either $f(x) = x$, or $f(x) = -x.$ \\~\\~\\
\textbf{11.}	Determine all functions $f: \Bbb{R}^+\to \Bbb{R}^+$ such that 
\begin{center}
$f(x+y+f(y)) = f(x) + 2y.$\\~\\
\end{center}
\textit{Solution:} First, we prove that the function is injective. Indeed, suppose that $f(a) = f(b).$ By considering $P(x,a)- P(x,b),$ we obtain $f(x+a+f(a)) = f(x+b+f(b)) + 2(a-b),$ or, defining $d = a-b,$ 
$f(x+d) = f(x) + 2d,$ for all $x > b+f(b).$ Considering $P(x, y+d),$ it follows that $d=0,$ i.e. $f(a) = f(b)$ implies $a=b,$ thus the function is indeed injective.
\\ Another way to obtain this is by considering $P(a,b)$ and then $P(b,a)$: \\
$f(a+b+f(b)) = f(a) + 2b$ \hfill (*)  \\ 
$f(a+b+f(a)) = f(b) + 2a$ \hfill (**) \\ 
By (*),(**)  (since $f(a)=f(b)$) we obtain $a=b.$ \\ 
Furthermore, by considering P$(x+f(x), y), P(y+f(y), x),$ we obtain:
$f(x+f(x)) = 2x+c$ for some constant $c,$ for all $x.$ \hfill (1)
\\ Considering $\displaystyle P\left(\frac{x}{2} + f\left(\frac{x}{2}\right), \frac{x}{2}\right)$ yields $\displaystyle f\left(x+2f\left(\frac{x}{2}\right)\right) = f\left(\frac{x}{2} + f\left(\frac{x}{2}\right)\right) + x = 2x + c,$ by (1).
\\Since f is injective, this means that $\displaystyle x+f(x) = x + 2f\left(\frac{x}{2}\right),$ hence $f(2x) = 2f(x).$ \hfill (2)
\\Considering $P(x+f(x), x),$ we have $2f(x+f(x)) = 4x+c;$ by (1), the LHS is also equal to $4x+2c.$ It follows that $c=0,$ i.e. $f(x+f(x)) = 2x.$ \hfill (3)\\~\\
\textit{Continuation 1:} Now we iterate (3) by taking $x\to f(x)$: this gives $f(f(x) + f(f(x)) = 2f(x) = f(2x),$ thus by injectivity we must have $f(x) + f(f(x)) = 2x$ (4).
The rest is just routine work on recurrence relations: Define $a_n = f^{(n)}(x)$ as the n-th iteration of $f$; then (4) gives $a_{n+2} = -a_{n+1} + 2a_n,$ which has characteristic equation $\lambda^2+\lambda -2 = 0,$ with roots $1$ and $-2.$ Hence $a_n = c_1\cdot 1 + c_2\cdot (-2)^n.$
Since all terms of the sequence are values of $f,$ they must be positive. However, if $c_2 \neq 0,$ then taking $n = 2k+1$ to infinity, or $n = 2k$ to infinity (depending on the sign of $c_2,$ in order to take $c_2\cdot(-2)^n$ to $- \infty),$ $a_n$ would eventually become negative, which is absurd; it follows that $c_2 = 0.$
Thus $a_n = c$ for all $n,$ and in particular $a_1 = c = a_0$ gives $f(x) = x.$ \\~\\
\textit{Continuation 2:} Consider $P(x+f(x), y): f(x+f(x) + y+f(y)) = 2(x+y),$ and call it $Q(x,y).$
\\ $\displaystyle Q(x,y), \ Q\left(\frac{x+y}{2}, \frac{x+y}{2}\right)$ give $\displaystyle f(x+f(x) + y+f(y)) = f\left(x+y + 2f\left(\frac{x+y}{2}\right)\right),$ and since $f$ is injective, this means that $\displaystyle \frac{f(x)+f(y)}{2}=f\left(\frac{x+y}{2}\right).$ This relation comes up often, so we prove the following lemma:\\~\\~\\
\textbf{Lemma (Jensen on the positive reals)}: Suppose that $f:\Bbb{R}^+\to\Bbb{R}^+$ is a function that satisfies
$\displaystyle P(x,y): \frac{f(x)+f(y)}{2}=f\left(\frac{x+y}{2}\right).$ Then $f(x) = ax+b,$ for nonnegative $a,b,$ not both $0.$ \\ \\
\textit{Proof:} We write $P(x,y)$ in the form $f(2x) + f(2y) = 2f(x+y).$ Fix some positive real $a.$\\ 
$P(x,c)$ gives $f(2x) = 2f(x+c) - f(2c),$ thus $P(x,y)$ becomes \\ $Q(x,y): f(x+c) + f(y+c) = f(x+y) + f(2c).$ \\ By the Cauchy-type lemma, or by considering $Q(x+c,y) - Q(x,y+c),$ it follows that \\ $f(x+c) = f(x) + d$ for some $d,$ for all $x>e.$ \\ Now, taking $Q(x,y)$ for $x,y > e,$ we have $f(x+y) = f(x) + f(y) + 2d - f(2c),$ which means that $g(x) = f(x) - (2d-f(2c))$ is additive for all sufficiently large $x,y;$ since it is also lower-bounded, we must have $g(x) = ax$ for all sufficiently large $x,$ that is $f(x) = ax+b$ for all $x>M.$
\\ Then $P(x,y)$ for $x>M$ means that $f(2y) = 2ay + b$ for all $y,$ hence $f(x) = ax+b$ for all $x. \square$ \\ \\
Returning to the problem, it follows that $f(x) = ax+b.$ Plugging this into the initial relation gives $f(x) = x$ as the only solution. \\~\\~\\
This problem has given us the opportunity to showcase two additional powerful weapons in our arsenal: recurrence relations and Jensen's functional equation. The following two examples illustrate the power of these two ideas assisted by our method.\\~\\~\\
\textbf{12.}	Determine all functions $f: \Bbb{R}^+ \to \Bbb{R}^+,$ such that 
\begin{center}
$f(x+f(y)) + f(y+f(x)) = 2(x+y).$ \\~\\
\end{center}
\textit{Solution:} $P(x,x)$ gives $f(x+f(x)) = 2x,$ or $f(g(x)) = 2x$ where $g(x) = x+f(x).$ Clearly, this shows that $f$ is surjective onto $\Bbb{R}^+.$ Furthermore, we have $g(g(x)) = f(g(x)) + g(x) = 2x+g(x).$ Thus, if $a_n = g^{(n)}(x),$ we obtain $a_{n+2} = a_{n+1}+2a_n,$ which has the associated characteristic equation $\lambda^2-\lambda-2 = 0,$ giving the solutions $2$ and $-1,$ hence $a_n = A\cdot 2^n + B\cdot(-1)^n.$ We immediately see that $A$ cannot be negative, since that would make an negative for sufficiently large $n,$ a contradiction; it follows that $A\geq 0.$ \\ 
At this point, we are stuck - we can't get anywhere else with this relation. However, if $n$ could also assume negative values (which would be allowed if $g$ were bijective) we would be done: for, if $B \neq 0,$ we could take $n$ even or odd (so that $B\cdot(-1)^n$ is negative), negative and sufficiently large (in magnitude) so that $A\cdot2^n$ is made arbitrarily small; again, this would make an negative, absurd; thus, we would have $a_n = A\cdot2^n = 2a_{n-1},$ implying $x+f(x) = g(x) = 2x,$ or $f(x)=x.$ So this motivates us to prove that $g$ is both injective and surjective. \\ 
First, $g(a) = g(b)$ would imply $2a = f(g(a)) = f(g(b)) = 2b,$ implying $a=b;$ it follows that $g$ is injective.
Also, suppose that $f(a) = f(b)$. Then $P(x,a) - P(x,b)$ gives $f(f(x)+a) - f(f(x)+b) = 2(a-b).$ Taking $a-b = d$ and remembering that $f(x)$ can assume any positive real value, we see that $f(x+d) = f(x) + 2d$ for all sufficiently large $x.$ $P(x+d, y) - P(x,y)$ gives $d=0,$ thus $f$ is indeed injective.
\\ Hence $f\displaystyle \left(g\left(\frac{f(x)}{2}\right)\right) = f(x),$ which means $\displaystyle g\left(\frac{f(x)}{2}\right) = x$ by injectivity; thus, $g$ is also surjective.
To conclude, $g$ is bijective as desired, proving that $f(x)=x$ is the only solution. \\ \\  \\
Actually, this idea we used for $f(g(x))$ in order to prove that $g$ is bijective can easily be generalized:\\~\\
\textbf{Composition Lemma}: Let $f(x)$ be an injective function and $h(x)$ a surjective function. If $g(x)$ is such that $f(g(x)) = h(x),$ then $g(x)$ is also surjective. The same holds if ''surjective'' is replaced with ''injective'' or ''bijective''.\\ 
\\ \textit{Proof:} If $h$ is injective, suppose that $g(a) = g(b).$ Then $h(a) = f(g(a)) = f(g(b)) = h(b)$, thus $a=b,$ meaning that $g$ is also injective.
\\If $h$ is surjective, then for every $x$ there exists some $s(x),$ such that $h(s(x)) = f(x).$ It follows that $f(g(s(x))) = f(x),$ and since $f$ is injective, this means $g(s(x)) = x,$ meaning that $g$ is also surjective.
\\If $h$ is bijective, we simply combine the previous arguments.\\~\\~\\
\textbf{13.} Determine all functions $f: \Bbb{R}^+ \to \Bbb{R}^+$ such that 
\begin{center}
$f(x+f(x)+2y) = 2x + 2f(f(y)).$ \\~\\
\end{center}
\textit{Solution:} First, by fixing $y$ and letting $x$ span the positive reals, we see that $f$ assumes all values on an interval $(2f(f(c)), +\infty).$ \\ 
We prove that $f$ is injective: indeed, assume that $f(a)=f(b)$ for $a>b,$ and let $d=a-b.$ \\ Then $P\displaystyle \left(a,\frac{x}{2}\right) - P\left(b,\frac{x}{2}\right)$ gives $f(x+a+f(a)) - f(x+b+f(b)) = 2d,$ or $f(x+d) = f(x) + 2d,$ for all $x > b+f(b).$ Now $P(x+d, y) - P(x,y)$ gives $d=0,$ absurd; it follows that f is injective. \\ \\
\textit{Continuation 1:} In order to get things equal, we consider $\displaystyle P\left(x, \frac{y+f(y)}{2}\right),$ which gives:\\
$\displaystyle f(x+f(x)+y+f(y)) = 2x+2f\left(f\left(\frac{y+f(y)}{2}\right)\right)$ \hfill (1) \\ 
By symmetry, the right-hand side is also equal to $\displaystyle 2y+2f\left(f\left(\frac{x+f(x)}{2}\right)\right).$ \\ It follows that $\displaystyle 2f\left(f\left(\frac{x+f(x)}{2}\right)\right) = 2x+c$ for some constant $c,$ for all $x.$ Thus (1) becomes \\ $f(x+f(x)+y+f(y)) = 2x+2y+c;$ call this $Q(x,y).$ \\
$Q(x,x)$ gives $f(2x+2f(x)) = 4x+c$. By the composition lemma we proved for problem 12, $Q(x,x)$ means that the expression $2(x+f(x))$ is injective and surjective onto an interval of the form $(M, +\infty).$
\\$ Q(x,y)$ can now be rearranged into:
\\$f(2x+2f(x)) + f(2y+2f(y)) = 2f(x+f(x)+y+f(y)).$ But since $x+f(x)$ can assume any sufficiently large value, we may write $x+f(x) = a, y+f(y) = b$ for arbitrary, sufficiently large $a, b$ to obtain $\displaystyle \frac{f(a)+f(b)}{2}=f\left(\frac{a+b}{2}\right)$ for all sufficiently large $a, b$. \\ By Jensen's functional equation, it follows that $f(x) = cx+d$ for all sufficiently large $x$ ($x>M$).
Take $y$ to be sufficiently large in the original equation (so that $\min\{y, cy+d\} > M).$ This yields $c(x+f(x)+2y) + d = 2x + 2c(cy+d) + 2d,$ which immediately gives $c=1, d=0,$ leaving $f(x)=x,$ for all $x.$ \\ \\
\textit{Continuation 2:} Let $a > b$ be arbitrary and set $d = a-b, k = f(a)-f(b).$ \\ Then $\displaystyle P\left(a,\frac{x}{2}\right) - P\left(b,\frac{x}{2}\right)$ yields $f(x+a+f(a)) - f(x+b+f(b)) = 2d,$ thus \\ $f(x+k+d) = f(x)+2d,$ for all $x > \min \{a+f(a), b+f(b)\}.$ \hfill (2) \\ Now $P(x,y+k+d) - P(x,y)$ gives $2d = f(f(y)+2d) - f(f(y)),$ from which it follows that $f(x+2d) = f(x)+2d$ for all sufficiently large $x$ (since $f(y)$ can assume all sufficiently large values).
\\Hence $P(x+d, y)$ writes $f(x+d+f(x+d)+2y) = 2x+2f(f(y))+2d = f(x+f(x)+2y+2d),$ and since f is injective, this means that the two arguments are equal, which implies $f(x+d) = f(x)+d.$
\\Hence (2) becomes $f(x+k) = f(x)+d = f(x+d)$ for all sufficiently large $x,$ and since $f$ is injective we see that $k=d,$ meaning that $a-b = f(a)-f(b).$ Since a,b were arbitrary, this means that $f(x) = x+c$ for some constant $c.$ Plugging this into the original equation yields $c=0,$ implying that $f(x) = x$ is the unique solution.\\~\\~\\
Actually, this last idea of considering variable $c$ and $d$ (usually certain differences between arbitrary numbers and the values of the function) can yield important equalities once we find the relation between $c$ and $d$ like we did in this problem. We see a few more examples to illustrate this idea:\\~\\~\\ 
\textbf{14.}	Determine all functions $f: \Bbb{R}^+ \to \Bbb{R}^+$  such that 
\begin{center}
$f(x+f(x)+2f(y)) = 2x+y+f(y).$\\~\\
\end{center}
\textit{Solution:} By fixing $y$ and letting $x$ vary, we can see that $f$ is surjective onto an interval of the form $[M, +\infty) .$ Fix some $a, b$ and let $k=f(a)-f(b)\geq 0, \ m = a-b.$ If $f(y_1) = f(y_2)$, $P(x,y_1), \ P(x,y_2)$ give $y_1 = y_2,$ thus the function is injective.
\\$P(a,y) - P(b,y): f(a+f(a)+x) - f(b+f(b)+x) = 2m,$ where $x = 2f(y)$ can be any number $\geq2M.$ That is, $f(x+k+m) = f(x) + 2m$ for any $x\geq  2M+b+f(b).$ \hfill (1)
\\$P(x+k+m, y) - P(x,y)$ gives $f(x+2m) = f(x) + 2k,$ for all sufficiently large $x$ (e.g. $x\geq 4M$)
\\$P(x+m, y): f(x+m+f(x+m)+2f(y)) = 2x+2m+y+f(y) = f(x+f(x)+2f(y)+k+m)$ by (1), for $x,y > 4M.$
\\Since $f$ is injective, this means that $x+m+f(x+m)+2f(y) = x+f(x)+2f(y)+k+m,$ or $f(x+m) = f(x)+k$ (2), for all $x > 4M.$
\\Now (1) becomes $f(x+k) = f(x) + 2m-k.$ Hence $P(x,y+k) - P(x,y)$ gives: $4k-2(2m-k) = 2m,$ thus $k=m.$ This means that $f(a) = a+C$ constant $C$ and all positive real $a.$ Replacing $f(x) = x+C$ back into the original equation, we obtain $f(x) = x$ as the unique solution. \\~\\~\\
\textbf{15.}	Determine all functions $f: \Bbb{R}\to \Bbb{R}$ such that 
\begin{center}
$f(x+y+f(x)^2) = x+xf(x)+f(y).$ \\~\\
\end{center}
$P(0,y)$ gives $f(y+f(0)^2) = f(y),$ thus $f(x)=f(x+T)$ for $T=f(0)^2.$ Taking $P(x+T,y),$ we obtain $Tf(x)=0$ for all $x,$ hence $f(0)=0.$ Also, take a number $u,$ such that $f(u)=0.$ \\ $P(u,x)$ gives $f(x+u) = f(x) + u$, whereas $P(x,u)$ now implies that $f(x+ f(x)^2) = f(x+ u +f(x)^2) = f(x+f(x)^2) + u,$ which gives $u=0.$ It follows that $f(x)=0 \iff x=0.$ \\ \\ 
\textit{Solution 1:} Denote $A = x+f(x)^2, B = x+xf(x).$
\\$P(x,0)$ gives $f(A) = B$, whereas $P(x, -x-f(x)^2)$ gives $f(-A) = -B = -f(A).$\\ By considering $P(A, -A)$ and $P(-A, A)$, we get $A + Af(A) -f(A) = -A + Af(A) + f(A),$ thus $A = f(A) = B,$ which gives $f(x)^2 = xf(x).$ Since $f(x)=0 \iff x=0,$ it follows that $f(x)=x$ for all real $x,$ which satisfies the initial relation.
\\ \\ Although this is a short solution, it is hard to motivate the use of $A$ and $B$ in it; instead, we can use the same idea on variable $c, d$ to work out an easy bash: \\  \\
\textit{Solution 2:} Fix $t$, and let $c = t+f(t)^2, \ d = t+tf(t)$. The initial relation then gives $f(x+c) = f(x)+d$ for all $x.$
\\$P(x+c, y) - P(x,y)$ gives $f(x+y+f(x)^2 + d^2 + 2df(x)) + d - f(x+y+f^2(x)) = xd + cf(x) + cd + c-d.$ 
\\Observe that $x+y+f(x)^2$ spans all of $\Bbb{R}$ as $y$ spans all of $\Bbb{R}$; thus we may set $x+y+f(x)^2 = Y$ equal to any real number, in order to get:
\\$Q(x,y): f(y+d^2+2df(x)) = f(y) +dx + cf(x) + cd + c-d $
\\$Q(0,y)$ gives $f(y+d^2) = f(y) + cd + c-d$, hence $Q(x,y)$ becomes \\ $R(x,y): f(y+2df(x)) = f(y) + xd + cf(x).$
\\Taking $x\to x+c$ in this one, we obtain $f(y+2df(x)+2d^2) = f(y) + xd + cf(x) + 2cd.$ But since $f(y+2d^2) = f(y) + 2(cd+c-d),$ we must have $2(cd+c-d) = 2cd$, which forces $c=d,$ i.e. $f(t)^2 = tf(t),$ which again gives $f(x) = x$ as the only solution.\\~\\~\\
\textbf{16.}	Determine all functions $f: \Bbb{R}^+\to \Bbb{R}^+$ such that 
\begin{center}
$f(x+f(x)+2y) = f(2x) + y + f(y).$\\~\\
\end{center}
\textit{Solution:}  We try to make arguments equal: if $x_0 > f(x_0)$ for some $x_0,$ then setting \\ $\displaystyle y\to \frac{x_0-f(x_0)}{2}$ would give $y+f(y) = 0,$ absurd. Thus $f(x)\geq x$ for all $x.$
\\Suppose that $f(a) = f(b)$ for some $a \geq  b.$
\\$\displaystyle P\left(\frac{a}{2}, x\right), P\left(\frac{b}{2}, x\right)$ gives $\displaystyle f\left(x + \frac{a}{2} + f\left(\frac{a}{2}\right)\right) = f\left(x + \frac{b}{2} + f\left(\frac{b}{2}\right)\right).$ \\ If $\displaystyle \frac{a}{2} + f\left(\frac{a}{2}\right) \neq \frac{b}{2} + f\left(\frac{b}{2}\right),$ denote their difference by $T.$ Then $f(x) = f(x+T)$ for all sufficiently large $x.$ Applying this to $P(x,y+T)$ for sufficiently large $y$ yields $T=0,$ so indeed $\displaystyle \frac{a}{2} + f\left(\frac{a}{2}\right) = \frac{b}{2} + f\left(\frac{b}{2}\right).$
\\Then  $\displaystyle P\left(a, \frac{b}{2}\right), P\left(b, \frac{a}{2}\right)$ give $f(2a) = f(2b).$ Considering $\displaystyle P\left(a, \frac{x}{2}\right), P\left(b, \frac{x}{2}\right)$ we get $f(x) = f(x+d)$ for all sufficiently large $x,$ where $d = a - b.$ Again, applying this to $P(x,y+d)$ for sufficiently large $y$ gives $d=0,$ thus the function is injective.
\\In general, $\displaystyle P\left(a, \frac{x}{2}\right) - P\left(b, \frac{x}{2}\right)$ gives $f(x+d) = f(x) + k$ for all $x > b+f(b),$ where \\ $d = a+f(a)-b-f(b) >0,\ k = f(2a) - f(2b).$
\\$P(x, y+d) - P(x,y)$ gives $2k = d+k$, hence $k=d.$ This means that $f(2a) = a+f(a) + t$ for some constant $t,$ for all $a.$ Also $f(x+k) = f(x)+k$ for all $x>M.$
\\Observe that $f(x+nk) = f(x) + nk$, and in particular $f(2nk) = f(nk) + nk$ (for $n$ such that $nk>M$). Also $f(2nk) = nk+f(nk)+t$ ; it follows that $t = 0.$  Thus $f(2x) = x+f(x)$ for all $x.$ 
\\The initial relation can now be written as $f(x+f(x)+2y) = f(2x) + y+f(y) = f(2x) + f(2y) = f(y+f(y)+2x)$; by injectivity, $x+f(x)+2y = y+f(y)+2x$, or $f(x) = x+c$ for all x; this satisfies for any nonnegative constant $c.$ \\~\\~\\
\textbf{17.} (Iran TST)	Determine all functions $f: \Bbb{R}^+ \to\Bbb{R}^+$ such that 
\begin{center}
$f(y)f(x+f(y)) = f(x)f(xy).$ \\~\\ 
\end{center}
\textit{Solution:} $P(1,y)$ gives $f(f(y)+1) = f(1).$ This motivates us to examine what happens if $f(a)=f(b)$: then $P(x,a), P(x,b)$ give $f(x) = f(kx)$ for all $x,$ where $\displaystyle k = \frac{b}{a} > 1. $
\\$P(kx,y)$ then gives $f(kx+f(y)) = f(x+f(y)).$ Now take a number $a,$ and define $\displaystyle t=\frac{ka+f(y)}{a+f(y)}.$
\\$P(x, a+f(y)),\ P(x, ka+f(y))$ give $f(x) = f(tx)$ for all $x.$ Observe that $t$ can assume all the values in $(1,k)$ for $k>1$ (the former by taking $a\to 0,$ the latter by taking $a\to  \infty$). But also $f(x) = f(x)$ (obviously) and $f(x) = f(kx), $ thus $f(x) = f(xy)$ for any $1\leq y \leq  k$ for some $k$ such that $f(x) = f(kx).$ But observe that $f(x) = f(kx)$ also implies that $f(k^2x) = f(kx) = f(x),$ and inductively $f(k^nx) = f(x),$ so $y$ can in fact assume any value (by choosing sufficiently large $n,$ so that $y<k^n).$ This means $f(x) = f(xy)$ for all $y\geq 1$ and all $x$; by symmetry, $f$ is constant. Conversely, any positive real constant function satisfies the given relations.\\~\\~\\ 
Although the next idea does not strictly belong to the cde method, it can certainly be motivated through it. We include it in this paper because of its beauty and usefulness:
If we have a function defined on the integers or a subset thereof and can inductively obtain a relation of the form $T(x,y,k): f(x+kA(y)) = \text{something},$ where $A(y)$ is a function of $y,$ then it is useful to consider $T(x,y,A(z))$ and $T(x,z,A(y)),$ which evaluate $f(x+A(y)A(z))$ in two ways.\\~\\~\\
\textbf{18.} (Italian IMO Camp)	Determine all integers $n,$ such that there exists a function $f: \Bbb{Z}\to \Bbb{Z}$ satisfying 
\begin{center}
$f(x+y+f(y)) = f(x)+ny.$ \\~\\
\end{center}
\textit{Solution:}  First of all, observe that $f(x) = 0$ satisfies this relation for $n=0.$ Assume now that $n \neq 0.$ \\ 
Then $P(x+y+f(y), y): f(x+2(y+f(y)) = f(x+y+f(y)) + ny = f(x) + 2ny.$ Inductively, we obtain the relation $T(x,y,k): f(x+k(y+f(y)) = f(x) + kny.$ \\ 
$T(x,y,z+f(z)), T(x,z,y+f(y))$ give $(z+f(z))y = (y+f(y))z$ for all $y,z$ (since $n \neq 0$). This means $x+f(x) = ax$ for some constant $a,$ for all $x.$ Taking $(y,z)=1,$ it follows that $y$ divides $f(y),$ implying that $a$ is an integer. This means $f(x) = cx$ for some integer $c,$ for all $x.$ Plugging this into the initial equation, $n$ must be of the form $c^2+c$ for some integer c.\\~\\~\\
\textbf{19.}	Determine all functions $f: \Bbb{Z}\to \Bbb{Z} $ such that 
\begin{center}
$f(2x + f(y)^2) = 2f(x) + yf(y).$\\~\\
\end{center}
\textit{Solution:}  First we make the RHS $0$: Taking $P(-2,-2)$, we have $f(u) = 0$, where $u = f(-2)^2 - 4. $
\\$P(x,u)$ gives $f(2x) = 2f(x)$, thus in particular $f(0)=0.$ Now we may create a coefficient imbalance: considering $P(2x, 2y)$, we get \\ $Q(x,y): f(x+f(y)^2) = f(x) + yf(y)$. Taking $x \to x+f(y)^2$ in this one gives $f(x+2f(y)^2) = f(x) + 2yf(y),$ and inductively we obtain $P(x,y,k): f(x+kf(y)^2) = f(x) + kyf(y).$
\\Now $P(x,y,f(z)^2), P(x,z,f(y)^2)$ give $f(y)f(z)[yf(z) - zf(y)] = 0$ for all y,z.
\\If $f(y)f(z) \neq 0$, then $yf(z) = zf(y)$; the second relation means $f(x) = ax$ for some nonzero constant $a$, and all $x$ such that $f(x) \neq 0.$
\\By $Q(0,y)$ we get $f(f(y)^2) = yf(y)$. Hence if $f(t) \neq 0$ for some $t,$ then also $f(f(t)^2) \neq 0$. Thus, $f(y) = at$ and $f(f(y)^2) = a(at)^2 = tf(t) = at^2$, which implies $a = 1$ or $a = -1.$ 
\\Note that $Q(x-f(y)^2, y)$ gives $f(x-f(y)^2) = f(x) - yf(y)$, so $f$ being a solution implies that $-f$ is also one. Therefore, we may assume that $a = 1.$ That is, for all $x,$ either $f(x) = 0$ or $f(x) = x.$
\\Let $A$ denote the set of zeroes of $f$, and let $B$ denote the set of nonzero fixed points of $f$. By the previous observation, their union is $\Bbb{Z}$. If $B$ is empty, then $f(x) = 0$ which clearly satisfies. Otherwise, consider $Q(x,b)$ for $b\in B$: $f(x+b^2) = f(x)+b^2.$
\\Hence $Q(-f(y)^2, y+b^2) - Q(-f(y)^2,y)$ gives $f(b^4+2b^2f(y)) = b^2(y+f(y)) + b^4$.
\\The LHS cannot equal $0$ (since we would need $b=0$, absurd by definition).
\\Thus $b^2(y+f(y)) + b^4 = b^4 + 2b^2f(y)$, which implies $f(y) = y$ for all $y.$ 
\\In sum, either $f(x) = 0$ for all $x,$ or $f(x) = x$ for all $x.$\\~\\~\\
\textbf{20.}	Determine all functions $f: \Bbb{Z}\to \Bbb{Z}$ such that 
\begin{center}
$f(x-y+f(y)) = f(x) + f(y).$\\~\\
\end{center}
\textit{Solution:}  Observe that $f(x) = 0$ everywhere is a solution; assume now that $f$ is not identically zero.
$P(x,x): f(f(x)) = 2f(x)$, thus $t\in \text{Im}f$ implies $2t\in \text{Im}f$; since we assumed that there is a nonzero value $t\in \text{Im}f,$ and inductively $2^kt \in \text{Im}f,$ it follows that $f$ has infinite image. \\ \\
\textit{Continuation 1:} Consider $P(x-(y-f(y)), y): f(x-2(y-f(y)) = f(x-(y-f(y)) + f(y) = f(x)+2f(y).$ Inductively, we get $T(x,y,k): f(x-k(y-f(y)) = f(x)+kf(y).$
\\$T(x,y,z-f(z)), \ T(x,z,y-f(y))$ give $(z-f(z))f(y) = (y-f(y))f(z),$ which implies $f(x)=ax$ for some constant $a$ and all $x.$ Plugging this into the original equation, we get $a=2$ or $a=0.$\\ \\
\textit{Continuation 2:} Assume now that $f$ is not injective; there must exist $a \neq b$ such that $f(a)=f(b)$. Then $P(x+b-f(b),a) - P(x+b-f(b),b)$ gives $f(x+T) = f(x)$ for all integers $x,$ where $T = a-b \neq 0.$ This would imply that $f$ is periodic, and since it is defined on the integers, it should have a finite image, clearly absurd by the previous considerations. Hence $f$ is injective.
Then $P(x,y), P(y,x)$ give $x-y+f(y) = y-x+f(x),$ thus $f(x) = 2x+c$ for some constant $c,$ for all $x.$ Plugging this to the initial equation gives $f(x)=2x$ for all $x.$\\ \\
\textit{Continuation 3:} Suppose that $f(x)$ does not equal $2x$ everywhere, and let $u = f(y_0) - 2y_0 \neq 0.$ $P(u,y_0)$ gives $f(u)=0.$
$P(x+u, y)$ gives $f(x+u) = f(x).$ Like before, this is absurd, thus $f(x) = 2x$ everywhere.\\~\\~\\
As we just saw, the idea about periodicity we referred to in the introduction and used in the solution of the previous problem can be particularly useful when dealing with functions defined on the integers, since then the function being periodic implies that is has a finite image.\\~\\~\\
\textbf{21.} (ISL 2015 A3)	Determine all functions $f: \Bbb{Z}\to \Bbb{Z}$ such that 
\begin{center}
$f(x-f(y)) = f(f(x)) - f(y) - 1 .$\\~\\
\end{center}
\textit{Solution:} First, assume that $f$ is constant: $f(x) = c.$ Then $c = -1,$ which is a solution.
Assume now that $f$ is not identically $-1$. \\ Trying to make arguments equal, we take $y\to f(x).$
This yields $f(x-f(f(x)) = -1,$ which means that there exists some $A$ such that $f(A) = -1.$ 
\\$P(x, A)$ gives $f(x+1) = f(f(x))$ (1). We can easily guess that $f(x) = x+1$ is also a solution; it suffices to prove that $f$ is injective in order to show that this is the unique solution. 
\\Assume that $f(a) = f(b)$ for $a > b$; by (1), $f(a+1) = f(f(a)) = f(f(b)) = f(b+1)$, and inductively $f(a+n) = f(b+n),$ for all natural numbers $n$; it follows that $f(x+T) = f(x)$ for all $x\geq b,$ $T=a-b.$ That is, $f$ is eventually periodic.
\\ Consider the function $g: \Bbb{Z}_{\geq b}\to \Bbb{Z}$  such that $g(x) = f(x)$ for all $x\geq b.$ Since $g$ is periodic and defined on a subset of the integers, it follows that $g$ has a finite image; in particular, it has a minimum and a maximum. Assume that $g(s) = M$ is the maximum and $g(t) = m$ is the minimum. (It isn't really necessary to define $g$; we are simply formalizing the result that $f$ is eventually periodic in order to treat the maximum and the minimum)
\\Due to (1), $P(x,y)$ can be written as: $g(x-g(y)) = g(x+1) - g(y) - 1$ for all $x,y, x-g(y)\geq b.$
Then $P(s+M, s)$ gives $2M+1 = g(s+M+1) \leq  M,$ thus $M \leq  -1.$
On the other hand, $P(t+M, t)$ gives $2m+1 = g(t+m+1)\geq m,$ hence $m \geq -1.$
\\It follows that $m\geq -1\geq M\geq m,$ which implies $m=M=-1,$ i.e. $g(x) = -1.$ That is, $f(x) = -1$ for all $x\geq b.$ 
\\Since $f(f(x)) = f(x+1)$ by (1), $P(x,y)$ can be written as:
$f(x-f(y)) = f(x+1) - f(y) - 1.$ Taking $y \geq b,$ we get that $f(x-1) = f(x+1)$ for all x, that is $f(x) = f(x+2k)$ for any integer $k.$ Making $k$ sufficiently large, so that $x+2k \geq b$ which forces $f(x+2k) = -1,$ we get $f(x) = -1$ for all integers, a contradiction since we assumed $f$ to be non-constant.
\\Thus, $f$ is indeed injective, which implies $f(x) = x+1. $\\~\\~\\
\textbf{22.} Determine all monotone functions $f: \Bbb{Z}\to \Bbb{Z}$ such that 
\begin{center}
$f(x^{2015} + y^{2015}) = f(x)^{2015} + f(y)^{2015}.$ \\~\\
\end{center}
\textit{Solution:} $P(0,0)$ gives $f(0) = 0.$ $P(x,0)$ gives $f(x^{2015}) = f(x)^{2015}.$ \hfill (1) \\ Thus the initial relation is equivalent to $Q(x,y): f(x^{2015} + y^{2015}) = f(x^{2015}) + f(y^{2015}).$
\\$Q(x, -x): f(x^{2015}) + f(-x^{2015}) = f(x)^{2015} + f(-x)^{2015} = 0,$ giving $f(-x) = -f(x).$
\\Relation (1) gives $f(1)^{2015} = f(1)$, hence $f(1)\in \{-1,0,1\}.$ Observe that if $f$ is a solution, then $-f$ is also a solution; thus, we may assume that $f(1) = 1$ or $f(1) = 0.$
\\If $f(1) = 0,$ then $Q(1,1)$ gives $f(2) = 0.$ By (1), we have $\displaystyle f(2^{2015^k}) = 0.$ Since $f$ is monotone and $\displaystyle f(2^{2015^k}) = f(0) = 0,$ this means $f(x) = 0$ for all $\displaystyle x\in \left[0, 2^{2015^k}\right].$ Taking $k$ to infinity, we have $f(x) = 0$ for all nonnegative $x.$ Since $f$ is odd, this means $f(x) = 0$ for all $x,$ which is a solution.
\\If $f(1) = 1,$ $Q(1,1)$ gives $f(2) = 2$, and by (1) we have $\displaystyle f(2^{2015^k}) = 2^{2015^k}.$ If $f$ is injective, $f$ must be a permutation of $\displaystyle \left[0, 2^{2015^k}\right],$ and since it is monotone, it must be the identity (as $\displaystyle 2^{2015^k}$ is always a fixed point). Hence $f(x) = x$ for all x in this interval, and taking $k$ to infinity implies $f(x) = x$ for all nonnegative $x.$ Since $f$ is odd, this means $f(x) = x$ for all x.
\\ Assume now that $f$ is not injective, i.e. there exist $a < b$ satisfying $f(a) = f(b).$ Since $f$ is odd, we may assume that $b$ is positive and $>1.$  By (1), we have $f(a^{2015}) = f(b^{2015})$; since $f$ is monotone, this means $f(x) = f(a^{2015})$ for all $x\in [a^{2015},b^{2015}].$ In particular, since $b^{2015} - a^{2015} > 1,$ we must have $f(a^{2015} + 1) = f(a^{2015}).$ By $Q(a,1),$ this implies $f(1) = 0,$ which is absurd. Thus $f$ is indeed injective, and our solutions are $f(x) = x,$ $f(x) = -x$ and $f(x) = 0.$\\~\\~\\
\textbf{23.}	Determine all functions $f: \Bbb{R}^+\to \Bbb{R}^+$  such that 
\begin{center}
$f(x+y+f(y)) = 4032x - f(x) + f(2017y).$ \\~\\
\end{center}
\textit{Solution:} $P(2017y, y)$ gives $f(2018y + f(y)) = 4032\cdot2017y,$ hence $f$ is surjective onto $\Bbb{R}^+.$
\\$P(x,a)-P(x,b)$ gives $f(x+a+f(a))-f(x+b+f(b)) = f(2017a) - f(2017b),$ giving \\ $f(x+C) = f(x) + D$ for all $x > b+f(b), \ C = a+f(a) - b - f(b) > 0, \ D = f(2017a) - f(2017b).$
\\$P(x+C, y) - P(x,y)$ gives $D = 2016C$, thus $f(2017a) = 2016(a+f(a)) + t$ for some constant $t$ and all $a.$ \\This can be written as: $\displaystyle a + f(a)=  \frac{f(2017a)-t}{2016}$. Since $f$ is surjective, this means that the expression $a+f(a)$ also assumes all positive real values greater than $\displaystyle -\frac{t}{2016}.$
\\We may now rewrite the initial expression by taking $y+f(y)$ to equal any positive real number $\displaystyle z > -\frac{t}{2016}.$   (We will be writing $z > M,$ where $\displaystyle M = \max \left\{0, -\frac{t}{2016}\right\},$ \\ thus $f(2017y) = 2016(y+f(y)) + t = 2016z + t).$
\\$Q(x,z): f(x+z) = 4032x – f(x) + 2016z + t,$ for all $x\in \Bbb{R}^+$ and $z > M.$
\\$Q(x,z), Q(z,x)$ give $4032x - f(x) + 2016z = 4032z - f(z) + 2016x,$ for all $x,z > M.$
\\Therefore $f(x) = 2016x + A$ for all $x>M$ and some constant $A.$ 
\\Take x arbitrary, $y > M$; then $x+y+f(y) > M$, $2017y > M$,\\ hence $2016(x+2017y+A) = 4030x - f(x) + 2016\cdot2017y$, which means that $f(x) = 2016x - 2016A,$ for all $x.$
\\Since this also holds for $x>M,$ we must have $2017A = 0,$ or $A=0$. These give $f(x) = 2016x,$ for all $x$, which clearly satisfies. \\~\\~\\ 
\textbf{24.}	Let $Q(x)\in \Bbb{R}[x]$ be a polynomial of degree $\geq 2.$ Determine all functions $f: \Bbb{R}\to \Bbb{ R}$ satisfying 
\begin{center}
$f(x+Q(y)+f(y)) = f(x).$ \\~\\
\end{center}
\textit{Solution:} $f(y) = -Q(y)$ is an obvious solution; suppose there exists another, thus $f(y_0) + Q(y_0) = T \neq 0$ for some $y_0.$ This gives $f(x) = f(x+T)$ for all $x.$
\\$P(x - Q(y) - f(y), y+T), P(x - Q(y) - f(y), y): f(x+Q(y+T)-Q(y)) = f(x)$ for all x,y. Since $Q(y+T)-Q(y)$ is a polynomial of degree $\deg Q - 1 \geq  1,$ it is surjective onto an interval of the form $[M, +\infty )$, or one of the form $(-\infty, M].$ \\ That is, $f(x) = f(y)$ for all $|x-y| \geq  M.$ But assume that $M > x-y > 0$; then, for sufficiently large $n,$ $(x+nT)-y > 0$. \\ Hence $f(x) = f(x+nT) = f(y)$ for all $x,y,$ which gives the solution $f(x) = c,$ for all x. \\~\\~\\ 
\textbf{25.}	Let $P(x) \in \Bbb{R}[x]$ be a polynomial taking positive reals to positive reals, of degree $\geq 2.$ Determine all functions $f: \Bbb{R}^+\to \Bbb{R}^+,$ such that 
\begin{center}
$f(f(x)+P(x)+2y) = f(x) + P(x) + 2f(y).$\\~\\
\end{center}
\textit{Solution:} Observe that $f(x) = c - P(x)$ is not a solution. Thus there exist $a,b$ such that $f(a)+P(a) > f(b)+P(b).$ Considering $\displaystyle S\left(a, \frac{x}{2}\right) - S\left(b, \frac{x}{2}\right),$ we get $f(x+c) = f(x)+c$ for all $x>e,$ where $c = f(a)+P(a) - f(b) - P(b) >0.$ 
\\For $x >e$, taking $S(x+c, y) - S(x,y),$ we get \\ $f(f(x) + P(x+c) + 2y) = f(f(x) + P(x) + 2y) + P(x+c) - P(x).$
\\By taking $\displaystyle y \to \frac{1}{2} (y - f(x) - P(x)),$ we get $f(y + u) = f(y) + u$ for all $u = P(x+c) - P(x)$, $y > f(x) + P(x)$. \\ Observe that $P(x+c) - P(x)$ is a polynomial of degree $\deg P-1 \geq  1$ and \\ $\displaystyle \lim_{x\to \infty}⁡[P(x+c)-P(x)]= + \infty$, hence $u$ can be any number on an interval of the form $[M, +\infty).$
\\Therefore, we get $f(x+y) = f(x) + y$ (1) for all $y > M, \ x > g(y) > M$ for some function $g$.
\\Using (1), as well as $S(x,y)$ for $y>M$,  $x > g(2y)$, we get \\ $f(f(x)+P(x)) + 2y = f(x) + P(x) + 2f(y)$. Fixing this sufficiently large $x$, it follows that $f(y) = y+t$ for some constant $t$, for all $y>M.$
\\For $x>M$, $S(x,y)$ gives $\displaystyle f(y) = y + \frac{t}{2}$ for all $y.$ Since this also holds for $y>M$,  $t=0$ and $f(x)=x,$ for all $x.$ \\~\\~\\
\textbf{26.}	Determine all functions f$: \Bbb{R}\to \Bbb{R}$ such that 
\begin{center}
$f(x+yf(x)) = f(x) + xf(y).$\\~\\
\end{center}
\textit{Solution:} $P(0,y): f(yf(0)) = f(0),$ meaning that either $f(0) = 0$ or $\displaystyle y \to \frac{y}{f(0)}$ means that $f$ is constant, which gives $f(x) = 0$ for all $x.$ Looking at non-constant solutions, we get $f(0) = 0.$ Also, if some $u$ satisfies $f(u)=0$,  $P(u,y)$ gives $u=0$ since $f$ is not identically zero.
\\So now we make the RHS zero: $P(-1,-1)$ gives $f(-1-f(-1)) = 0,$ thus $f(-1) = -1.$
\\$P(-1,x): \ f(-x-1) = -f(x)-1.$  \hfill (1)
\\$P(x,-1): \ f(x-f(x)) = f(x)-x.$ \hfill 
\\Clearly, $f(x)=x$ for all $x$ is a solution. If there is another solution, then letting $x-f(x) = A$ be nonzero means that $ f(A) = -A$. Now $P(A,1)$ gives $0 = A(f(1)-1),$ so $f(1) = 1.$
\\$P(1,x): f(x+1) = f(x)+1$. Now (1) writes $f(-x) = -f(x).$ \\ Also, $f(n) = n,$ for all integers $n.$
\\$P(n,x)$ gives $f(nx) = nf(x)$ for all integers $n$ and reals $x.$ \\ Now we may create a coefficient imbalance:
\\$P(x,2y), P(x,y)$ give $f(x+2yf(x)) - f(x) = 2xf(y) = 2(f(x+yf(x)) - f(x)).$ \\ Thus $f(x+2yf(x)) + f(x) = 2f(x+yf(x)).$
\\But now observe that, since $x \neq 0$ implies $f(x) \neq 0,$ we may set $\displaystyle y \to \frac{y}{f(x)},$ for $x \neq 0,$ to obtain $f(x+2y) + f(x) = 2f(x+y) = f(2x+2y),$ for all $x \neq 0$; this holds identically for $x=0$ too.
\\But letting $x+2y = u, x = v$ where $u,v$ may be any real numbers, we get $f(u+v) = f(u)+f(v).$
\\Using this in the original equation, we obtain $f(yf(x)) = xf(y)$ for all real $x,y$. For $y \neq 0$, letting x span $\Bbb{R}$ means that $f$ is surjective. Taking $y=1$ gives $f(f(x)) = x,$ hence $f(yf(x)) = f(y)f(f(x))$ for $x,y,$ and since f$(x)$ is surjective we may write $f(x) = z$ for any real $z.$
\\So $f(yz) = f(y)f(z),$ and since $f(x+y) = f(x)+f(y),$ this means that $f(x) = x.$\\~\\~\\
\textbf{27.}	Determine all functions $f:\Bbb{R}^+\to \Bbb{R}^+$  such that 
\begin{center}
$f(xf(x)+f(y)) = f(x)^2 + y.$\\~\\
\end{center}
\textit{Solution:} The iteration lemma gives $f(y+K) = f(y)+K,$ where $K=xf(x)+f(x)^2.$ Inductively, we obtain $f(y+nK)= f(y) + nK$, thus $y>mK$ implies $f(y)>mK.$ \hfill (1)
\\ Now iterate with $y\to zf(z)+f(y)$ to obtain:
$f(y+a) = f(y)+b,$ where $a = xf(x)+f(z)^2,\ \ b=zf(z)+f(x)^2.$ Assume, for contradiction, that $a \neq b$; without loss of generality let $a>b$. Take a sufficiently large natural number $N$, such that $y+Na > f(y) + Nb + 3K,$ and for that $N$ let $\displaystyle m = \left\lfloor \frac{y+Na}{K}\right\rfloor - 1.$ This gives:
$y+Na > mK > f(y) + Nb = f(y+Na),$ yet $y+Na > mK$ should imply that $f(y+Na) > mK,$ by (1).
\\Combining these relations, we obtain: $mK > f(y+Na) > mK,$ clearly absurd; it follows that indeed $a=b$, which implies that $f(x)^2 = xf(x) + c$, where $c$ is some constant. \hfill (2) \\ The rest is easy: 
By substituting $x\to x+K$ in (2), and then using (1), we obtain:
\\$(x+K)(f(x)+K) + c = f(x+K)^2 = (f(x)+K)^2$, hence $K(f(x)+x) + c = 2Kf(x)$. But the initial relation shows that $f$ is unbounded, thus $K$ is unbounded; this means that $f(x)=x. $\\~\\~\\
\textbf{28.} Determine all functions $f: \Bbb{R}^+ \to \Bbb{R}^+$ such that 
\begin{center}
 $f(1+xf(y)) = yf(x+y).$ \\~\\
 \end{center} 
\textit{Solution:} This is a particularly hard problem. We give a solution in the cde spirit. \\
\\First of all, observe that $y>1$ implies $f(y) \leq 1,$ else $\displaystyle P\left(\frac{y-1}{f(y)-1},y\right)$ forces $y=1,$ absurd.
\\$1+xf(y) > 1$ implies $yf(x+y) = f(1+xf(y)) \leq 1,$ thus $f(x)\leq  \frac{1}{x-\varepsilon}$ for any positive $\varepsilon$. \\ Taking $\varepsilon \to 0,$ we get $\displaystyle f(x) \leq  \frac{1}{x},$ for all $x.$ \hfill (1) \\
\\Set $f(1) = a.$ $P(x,1)$ gives $f(1+ax) = f(x+1),$ yielding inductively $f(a^nx + 1) = f(x+1)$ (setting $\displaystyle x\to  \frac{x}{a^n}$, we see that this also holds for negative $n$) \\ It follows that $\displaystyle f(x+1) \leq \frac{1}{1+a^nx}$ for all $x,$ and taking $n\to \infty$  makes $f(x)$ arbitrarily small if $a > 1,$ contradiction. \\ Also taking $n\to -\infty$  gives a similar contradiction for $a < 1.$ Thus $a = f(1) = 1.$ \\
\\Also $\displaystyle P\left(\frac{x}{f(y)}, y\right)$ gives $\displaystyle f(x+1) = yf\left(y + \frac{x}{f(y)}\right) \leq \frac{y}{y+ \frac{x}{f(y)}}$ for all $x,y.$ \\ This can be written in the form $\displaystyle f(y) \geq \frac{b}{y},$ where $\displaystyle b = \frac{a}{\frac{1}{f(a+1)}-1}.$ \hfill (2)  \\
\\Now we prove that the function is injective. Suppose that $f(u) = f(v)$; then $P(x,u), \ P(x,v)$ give $\displaystyle f(x+u) = \frac{v}{u } f(x+v)$, thus $f(x+d) = cf(x)$ for $\displaystyle d = u-v >0, \ c =\frac{v}{u} <1,$ and all $x > v.$
\\Hence $f(x+nd) = c^nf(x).$ This implies $\displaystyle c^nf(x) = f(x+nd) \geq  \frac{b}{x+nd},$  by (2).
\\It follows that $\displaystyle f(x) \geq  \left(\frac{1}{c}\right)^n \frac{b}{x+nd}.$ \\ The RHS is $\displaystyle \mathcal{O}\left(\frac{e^n}{n}\right)$ which clearly goes to infinity for $n\to \infty $, a contradiction. Thus f is injective.\\
\\Now $P(1,x)$ gives $f(f(x)+1) = xf(x+1),$ whereas $\displaystyle P\left(1+f(x) - \frac{1}{x}, \frac{1}{x}\right)$ \\ (for $\displaystyle x>1,$ $\displaystyle 1+f(x)-\frac{1}{x}$ is positive) gives \\ $\displaystyle f\left(1+ f\left(\frac{1}{x}\right)\left(1+f(x) -\frac{1}{x}\right)\right) = \frac{1}{x}f(f(x)+1) = f(x+1).$ \\Since $f$ is injective, this means that $\displaystyle x = f\left(\frac{1}{x}\right) \left(1+f(x) - \frac{1}{x}\right).$ 
\\Since $\displaystyle f\left(\frac{1}{x}\right) \leq  x$ for all $x,$ this means $\displaystyle f(x)\geq \frac{1}{x},$ thus $\displaystyle \frac{1}{x}\geq f(x)\geq \frac{1}{x}.$ \\ It follows that $\displaystyle f(x) = \frac{1}{x},$ for all $x>1.$\\
\\From the initial equation we get $\displaystyle f(x+y) = \frac{1}{y+xyf(y)}$ for all $x,y.$ For $x>1,$ so that $x+y>1,$ we get $y+xyf(y) = x+y$ for all $y,$ hence $\displaystyle f(y) = \frac{1}{y}$ for all $y.$\\~\\

\textbf{3.	Conclusion} \\~\\
Having presented several examples that demonstrate how to apply the cde method, often in conjunction with other techniques, we should also observe that the takeaway from this article doesn't have to be a specific solution, but rather a strong insight into creating 'coefficient imbalances' that let us split up terms and make progress. As exemplified by the main point, taking the difference $P(x+c,y)- P(x,y)$, we are more interested in obtaining '+ something' than its actual value; it suffices to prove things 'for sufficiently large arguments'; in general, we first describe general forms and then refine them. The other important application of the technique can be described as proving injectivity before employing other tools; together, all these ideas account for a powerful arsenal, even against demanding problems. \\ 
Happy cde-ing! \\~\\~\\

\newpage

\textbf{4.	Practice problems} \\~\\
\textbf{1.} Determine all functions $f:\Bbb{R}^+\to \Bbb{R}^+$  such that $$f(x+f(y))=f(x)+2xy^2+y^2 f(y).$$ \\
\textbf{2.} Solve example problem 5 without using the cde method.   \\ \\ \\
\textbf{3.} Determine all functions $f:\Bbb{R}^+\to \Bbb{R}^+$  such that $$f(f(x)+2y)=f(2x+y)+2y.$$ \\ 
\textbf{4.} Determine all functions $f:\Bbb{R}^+\to \Bbb{R}^+$  such that $$f\left(\frac{x+f(x)}2+y\right)=2x-f(x)+f(f(y)).$$   \\
\textbf{5.} Determine all functions $f:\Bbb{R}^+\to \Bbb{R}^+$  such that $$ f(x+y+f(y))+f(x+z+f(z))=f(2f(z))+f(2y)+f(2x).$$    \\ 
\textbf{6.} Determine all functions $f:\Bbb{R}^+\to \Bbb{R}^+$  such that $$f(xy+f(x))=xf(y)+f(x).$$ \\ 
\textbf{7.} Determine all functions $f:\Bbb{R}^+\to \Bbb{R}^+$  such that $$f(x+yf(x))=f(x)+xf(y).$$ \\~\\~\\~\\~\\

\newpage

\textbf{5.	References} 
\begin{enumerate}
    \item Functional Equations A Problem Solving Approach, B.J. Venkatachala, Prism Publications 2008
    \item Functional Equations and How to Solve Them, C.G. Small, Problem Books in Mathematics 2007
    \item Functional Equations in Mathematical Olympiads (2017 - 2018): Problems and Solutions (Vol. I), A.H. Parvardi, 2018
    \item Functional Equations, Marko Radovanovic, The IMO Compendium Group, 2007
    \item Topics in Functional Equations, Andreescu et al., XYZ Press, 2012
    \item Introduction to Functional Equations, Evan Chen, 2016
    
\end{enumerate}

$$ $$ 
Athanasios Kontogeorgis		               \hfill		    Rafail Tsiamis          \\
Aristotle University of Thessaloniki       \hfill           Harvard College  \\
Thessaloniki, Greece		             
\hfill		    Boston, Massachusetts      \\
athkonto@gmail.com		                   \hfill		    rtsiamis@college.harvard.edu         \\

\end{document}